\documentclass[12pt,twoside]{amsart}
\usepackage{hyperref}
\usepackage{amsthm, amsmath, amscd, amssymb,centernot}
\usepackage[all]{xy}
\usepackage[T1]{fontenc}
\usepackage[left=2cm,top=2.5cm,bottom=3cm,right=3cm]{geometry}

\usepackage{amsthm}
\usepackage{tikz-cd}
\usepackage{amssymb}
\usepackage{enumerate}
\usepackage{amsmath} 
\usepackage[mathscr]{euscript}
\newtheorem{thm}{Theorem}[section]
\newtheorem{lem}[thm]{Lemma}

\newtheorem{xrem}[thm]{Remark}
\newtheorem{exm}[thm]{Example}
\newcommand{\reals}{\mathbb R}

\newcommand{\str}{\mathcal{O}}

\newcommand{\proj}{\mathbb{P}}

\theoremstyle{plain}
\numberwithin{equation}{section}

\newtheorem*{theorem*}{Theorem}

\theoremstyle{definition}

\begin{document}

\title{Weyl and Zariski chambers on projective surfaces}
\author[Krishna Hanumanthu]{Krishna Hanumanthu}
\address{Chennai Mathematical Institute, H1 SIPCOT IT Park, Siruseri, Kelambakkam 603103, India}
\email{krishna@cmi.ac.in}

\author[Nabanita Ray]{Nabanita Ray}
\address{School of Mathematics, Tata Institute of Fundamental Research, Homi Bhabha Road, Mumbai 400005, India}
\email{nray@math.tifr.res.in}

\thanks{First author was partially supported by a grant from Infosys
        Foundation and by DST SERB MATRICS grant MTR/2017/000243.}

\date{April 16, 2020}
\maketitle

\begin{abstract}
Let $X$ be a nonsingular complex projective surface. The Weyl and Zariski
chambers give two interesting decompositions of the big cone of $X$. 
Following the ideas of \cite{BF} and \cite{RS},  
we study these two decompositions and determine when a Weyl chamber is
contained in the interior of a Zariski chamber and vice versa. We also 
determine when
a Weyl chamber can intersect non-trivially with a Zariski chamber.  

 \end{abstract}

\vskip 4mm

\section{Introduction}

Let $X$ be a nonsingular projective surface. A divisor $D$ on $X$ is
called \textit{big} if the global sections of $mD$ determine a birational
morphism for $m \gg 0$ . The classes of all big divisors in the
N\'eron-Severi space $N^1(X)_{\reals}$ form a cone called the \textit{big cone} of $X$. 
The big cone is open and contains the ample cone of $X$. Studying 
the big cone of $X$ is useful in understanding the geometry of
$X$. The intersection of the big and nef cones is particularly interesting.  

A useful way to study the big cone is to 
decompose it into chambers and study individual chambers. 
One instance of such a decomposition is given by 
\textit{Zariski chambers}. These have been defined in \cite{BKS} as
the set of big divisors for which the negative part of the Zariski
decomposition is constant. Further, the Zariski chambers  are rational 
locally polyhedral subcones and the big cone admits a locally
finite decomposition into Zariski chambers. 
This decomposition turns out to have
significant geometric implications. 
On each of the Zariski chambers, the volume function of the divisors
is given by a single quadratic polynomial. Moreover, the stable base
loci are constant in the interior of each Zariski chamber. See
\cite{BKS} for more details.

It is thus an interesting problem to study the 
Zariski chamber decomposition of the big cone of a given surface. 
For example, one can ask for the number of chambers in the
decomposition. 
This question has been answered for Del Pezzo surfaces and some other surfaces with higher Picard number
in \cite{BFN} and \cite{BS}. The number of Zariski chamber is an
interesting geometric invariant of $X$. 

Another interesting decomposition of the big cone is given by
\textit{Weyl chambers}. 
Traditionally the (simple) Weyl chambers are defined on $X$ if the
only irreducible and reduced curves of negative self-intersection are (-2) curves. 
In this situation, \cite{BF} describes the the Weyl chamber decomposition for K3 surfaces.
But in \cite{RS}, the authors extend the definition of Weyl chambers
to arbitrary surfaces.

Since the Zariski and Weyl chambers are defined for any nonsingular
surface, it is interesting to study how they are related to each
other. The Zariski chambers are, in general, neither open nor closed,
but the Weyl chambers are always open. So it is natural to ask when
the interior of a Zariski chamber coincides with a Weyl chamber. 
This problem was studied for K3 surfaces in \cite{BF} and the
authors establish a 
necessary and sufficient condition for the interior of every Zariski
chamber to coincide with a Weyl chamber; see \cite[Theorem
2.2]{BF}. This result is generalized to an arbitrary nonsingular surface
in \cite[Theorem 3]{RS}.  

In this note, we continue the comparison of Zariski and Weyl chambers started by 
\cite{BF} and \cite{RS}. In Section \ref{prelims}, we 
prove that there is a bijection between the set of
Zariski chambers and the set of Weyl chambers on any nonsingular
projective surface (Theorem \ref{thm2.4}). This generalizes
\cite[Theorem 1.3]{BF}, which proved the bijection for K3 surfaces.

Our first main results (Theorems \ref{prop3.1} and
\ref{prop3.2})  give necessary and
sufficient conditions for a \textit{specific} Weyl chamber to be contained in a
Zariski chamber and for the interior of a \textit{specific} Zariski chamber to be
contained in a Weyl chamber.  We note that this result has been obtained
in \cite{BF} for K3 surfaces and our proofs are similar to the proofs of
\cite{BF}. We also show that our main results 
imply the main theorem of \cite{RS}. 
In another result (Theorem
\ref{thm3.7}), we give a necessary and sufficient condition for a
Zariski chamber and a Weyl chamber to have non-empty intersection. At
the end, we give some examples which illustrate our results.  

In Section \ref{prelims}, we recall the definitions of Weyl and
Zariski chambers and their basic properties before proving some new
results. In Section \ref{main-results}, we prove our main results
about comparing the Weyl and Zariski chambers and
give some examples.

Throughout this note,  we work over the field $\mathbb{C}$ of complex
numbers. If $X$ is a nonsingular projective surface, a
\textit{negative curve} on  $X$ is an irreducible and reduced curve $C$ satisfying
$C^2 < 0$.

\subsection*{ Acknowledgement} Part of this work was done when the
first author visited the Tata Institute of Fundamental Research,
Mumbai. He is grateful for the hospitality.

\section{Preliminaries}\label{prelims}

We start by recalling the main objects that we study. 

\subsection{Big cone}
A line bundle $L$ on an irreducible projective variety $X$ is
\textit{big} if the mapping $\phi_m: X\dashrightarrow
\mathbb{P}(H^0(X, L^{\otimes m}))$, defined by $L^{\otimes m}$ is
birational onto its image for $m\gg 0$. A Cartier divisor $D$ on $X$
is big if $\mathcal{O}_{X}(D)$ is big. Big divisors are precisely the
divisors whose volume is positive. Another characterization of big
divisors is given by the following: a divisor $D$ is big if and only
of $mD$ is numerically equivalent  to $A+E$ for an ample divisor $A$
and an effective divisor $E$ and some positive integer $m$.  In
particular, bigness is a property of the numerical class of $D$. 

The notion of bigness can be extended to $\reals$-divisors. An
$\mathbb{R}$-divisor $D\in \text{Div}_{\mathbb{R}}(X)$ 
is big if it can be written in the form $D =\sum a_iD_i$, where each $D_i$ is a big integral divisor and $a_i$ is a positive real
number.  Since bigness is preserved under numerical equivalence, 
we can talk about big classes in the real N\'{e}ron Severi group
$N^1(X)_{\mathbb{R}}$. The convex cone of all big $\mathbb{R}$-divisor
classes in $N^1(X)_{\mathbb{R}}$ is called the \textit{big cone} 
and it is denoted by $\text{Big}(X)\subset N^1(X)_{\mathbb{R}}$. 
For further details on big divisors and the big cone, see \cite{L}.

\subsection{Zariski Decomposition}
Zariski \cite{Z} showed that any effective divisor on a surface can be written
uniquely as a sum of a nef divisor and some negative curves. This has
been generalized to a larger class of divisors than effective divisors
and is known as the \textit{Zariski
  decomposition}; see \cite{F,B}. 

The closure of the big cone is called the \textit{pseudoeffective
  cone} and divisors in that cone are called pseudoeffective. Let $D$
be a pseudoeffective  $\mathbb{R}$-divisor on a nonsingular complex
projective surface $X$. Then there exist $\mathbb{R}$-divisors $P_D$ and $N_D$ such that
\begin{align*}
 D=P_D+ N_D,
\end{align*}
and the following conditions hold:
\begin{enumerate}
\item $P_D$ is nef,
\item either $N_D=0$ or $N_D=\sum\nolimits _{i=1}^{r} a_i E_i$, where
  $a_i > 0$ 
  and the intersection matrix \\$(E_i\cdot E_j)_{1 \le i,j\le r}$ is
  negative definite, and 
\item $P_D$ is orthogonal to each of the components of $N$, i.e.,
  $P_D\cdot E_i=0$, for $i=1,\cdots, r$.
\end{enumerate}
$P_D$ and $N_D$ are called the \textit{positive part} and the \textit{negative part} of $D$ respectively.
They are are uniquely determined by $D$.

It follows from the property (2) of the definition that the negative
part of the big divisor is either trivial 
or is supported by negative curves.

\subsection{Zariski chambers}
The variation of the Zariski decomposition over the big cone
$\text{Big}(X)$ leads to a partition of the big 
cone, given by subcones for which the negative part of the Zariski
decomposition is constant. Each such subcone is called a
\textit{Zariski chamber}. 
A formal definition of a Zariski chamber is given below. 

Let $D$ be a big divisor and let $D=P_D+N_D$ be the Zariski
decomposition of $D$. We first define two sets of curves associated to $D$
as follows: 
\begin{center}
 Null($D$) = $\{C \subset X \mid C$ is an irreducible curve with $D\cdot C=0\}$, and
\end{center}
\begin{center}
 Neg($D$) = $\{C \subset X \mid C$ an irreducible component of $N_D\}$.
\end{center}
Clearly, Neg$(D)\subseteq$ Null$(P_D)$. It may happen that 
Neg$(D)\ne $ Null$(P_D)$.

Suppose now that $P$ is a big and nef divisor.  The \textit{Zariski
  chamber} $\sum_P$ associated to $P$ is defined as

\centerline{$\sum_P=\Big\{D\in \text{Big}(X) \mid \text{Neg}(D) = \text{Null}(P)\Big\}.$}

The interior of $\sum_P$ is given by 
\begin{center}
 $\text{int} ({\sum_P}) = \{ D\in \text{Big}(X)\mid $ Neg($D$) =
 Null($P$) = Null($P_D$)$\}$.
\end{center}

The above result is proved in \cite[Proposition 1.8]{BKS}. It is also
known that the sets $\sum_P$ cover $\text{Big}(X)$; see \cite[Proposition 1.6]{BKS}.


If $H$ is an ample divisor, then the interior of the chamber $\sum_H$
is the ample cone and its closure is the nef cone. The Zariski chamber
$\sum_H$ for any ample divisor $H$ is called the \textit{nef
  chamber}.

We recall the following useful observation. 
\begin{lem}\cite[Lemma 1.1]{BFN}\label{lem2.1}
 The set of Zariski chambers on a smooth projective surface $X$ that
 are different  from the nef chamber is in bijective correspondence
 with the collection of sets of reduced
divisors on $X$ whose intersection matrix is negative definite.
\end{lem}
We define the following two sets. 
\[
\mathscr{I}(X) =\left\{
C \subset X\biggm| \begin{array}{l}
C \text{~is an irreducible and reduced}\\ 
\text{curve satisfying~} C^2 < 0
\end{array}
\right\},\text{~and}
\]
\[
\mathscr{Z}(X) =\left\{
S \subseteq \mathscr{I}(X) \biggm| \begin{array}{l}
S \text{~is finite and the intersection}\\ 
\text{matrix of~ } S \text{~is negative definite}
\end{array}
\right\}.
\]

The following remark is easy to prove and we use it repeatedly in our
arguments:  if $S \in \mathscr{Z}(X)$ and $S' \subset S$, then $S' \in
\mathscr{Z}(X)$.

To every $S \in \mathscr{Z}(X)$, we associate a set of big divisors as
follows: 
$$Z_S = \{D\in \text{Big}(X)\mid \text{Neg}(D)=S\}.$$
In words, $Z_S$ consists precisely of those big divisors whose
negative part is supported on $S$. 

Then $Z_S$ are precisely the Zariski chambers. Indeed, let  
$P$ be a big and nef divisor. We note that 
$\text{Null}(P) \in \mathscr{Z}(X)$: since Null$(P) = \text{Null}(mP)$, 
for any positive integer $m$, we may suppose that $P = A + E$, where $A$ is 
an ample divisor and $E$ is an effective divisor. If $C \in \text{Null}(P)$, 
then $N_E\cdot C < 0$. So $C$ is a component of $N_E$. Since the
irreducible components of $N_E$ have a negative definite intersection
matrix, so will the elements of $\text{Null}(P)$. 
Now let $S: = \text{Null}(P)$. Then we get
\begin{center}
 $\sum_P=\{D\in \text{Big}(X)\mid$ Neg$(D)=S\}=Z_S$.
\end{center}
Conversely, given $S \in \mathscr{Z}(X)$, it can be shown 
that there exists a big and nef divisor $P$ such that $\text{Null}(P) = S$, 
so that $\sum_P = Z_S$.

By Lemma \ref{lem2.1}, it follows that the set of Zariski chambers is in
bijective correspondence with the elements of $\mathscr{Z}(X) \cup
\emptyset$. Note that $Z_{\emptyset}$ is precisely the nef chamber in
the big cone of $X$. For $S, S' \in \mathscr{Z}(X)$, the Zariski chambers $Z_S$ and $Z_{S'}$ are either disjoint or identical. 

The following useful remark will be used in the proofs. 
\begin{xrem}\rm
Let $C_1,\cdots,C_r \in \mathscr{I}(X)$ and let $P$ be a nef divisor
such that $P+c_1C_1+\cdots+c_rC_r \in \text{Big}(X)$ for some positive
real numbers $c_1,\cdots,c_r$. Then $P+b_1C_1+\cdots+b_rC_r \in
\text{Big}(X)$ for any positive real numbers $b_1,\cdots,b_r$. We give
a brief argument for this fact below. 

First, note that $P, C_1,\cdots, C_r$ belong to the closure of the big cone,
since the closure $K$ of the big cone is the closure of the cone of
effective curves which contains the nef cone. Let $b_1,\cdots,b_r$
be any positive real numbers. If  $P+b_1C_1+\cdots+b_rC_r \notin
\text{Big}(X)$, then  $P+b_1C_1+\cdots+b_rC_r$ is contained in the
boundary of $K$. Note that the interior of $K$ is precisely
$\text{Big}(X)$ and it is disjoint with the boundary of $K$. If
$P+b_1C_1+\cdots+b_rC_r$  is contained in the boundary of $K$, then it
spans a ray $R$ in the boundary. This means that $P,C_1,\cdots, C_r
\in R$ which in turn means that $P+c_1C_1+\cdots+c_rC_r$ is contained
in the boundary of $K$. This contradicts the hypothesis that
$P+c_1C_1+\cdots+c_rC_r$ is a big divisor. 
\end{xrem}

For further details about the Zariski chambers, see \cite{BKS}.

\subsection{Weyl chambers}
Each curve $C\in \mathscr{I}(X)$ defines a hyperplane in the 
N\'{e}ron-Severi space $N^1(X)_{\mathbb{R}}$ of $X$ as follows:
 $$C^{\bot}=\{D \in N^1(X)_{\mathbb{R}} \mid D\cdot C=0\}.$$

These hyperplanes give a decomposition of the big cone
Big($X$). Namely, we consider the connected components of 
$$\text{Big}(X)\setminus \bigcup_{C\in \mathscr{I}(X)} C^{\bot}.$$
The connected components of the set $\text{Big}(X)\setminus \bigcup_{C\in
  \mathscr{I}(X)} C^{\bot}$ are 
called the simple \textit{Weyl chambers} of $X$.

Traditionally, the (simple) Weyl chambers are studied only if all the negative
curves on $X$ are $(-2)$-curves, i.e., smooth rational curves $C$
satisfying $C^2 = -2$; 
see \cite{BF}. But in \cite{RS}, 
authors studied Weyl chambers on arbitrary surfaces. 
In this note, we adopt the same approach. First we give a
convenient alternate
definition of Weyl chambers. This characterization was used in \cite{BF} 
for K3 surfaces.

Let $S\in \mathscr{Z}(X)\cup\{\emptyset\}$. Corresponding to $S$, we define
the set $W_S$ as follows: 
\begin{center}
 $W_S=\{ D\in \text{Big}(X) \mid D.C< 0$ $\forall C\in S$ and $D.C>0$ $\forall C\in \mathscr{I}(X)\setminus S \}$.
\end{center}
We show below that $W_S$ are precisely the Weyl chambers of $X$. 
Note that $W_{\{\emptyset\}}$ is the ample cone in Big($X$).
From the definition, it is clear that each $W_S$ is a convex set,
hence connected.

We first prove the following lemma. 
\begin{lem}\label{lem2.2}
The following properties hold:
\begin{enumerate}  
\item $W_S$ is an open set for every $S \in  \mathscr{Z}(X)\cup\{\emptyset\}$.
\item For $S, S' \in  \mathscr{Z}(X)\cup\{\emptyset\}$, $W_{S'} \cap W_S=\emptyset$ if and only if $S'\neq S$.
\end{enumerate}
\begin{proof}
(1)  To show $W_S$ is open, it is enough to show that for any $D\in
W_S$, and for any $E\in N^1(X)_{\mathbb{R}}$, $mD\pm E\in W_S$ for all
sufficiently large integers $m$. 
 
As the big cone is open, there exists an integer $m \gg 0$ such
that $mD\pm  E\in \text{Big}(X)$.

For $C \in S$, we need to ensure $(mD\pm E) \cdot C < 0$. Since
$D\cdot C < 0$ and the set of numbers $\{E\cdot C \mid C \in S\}$ is
finite, we can choose a bigger $m$ (if necessary) such that $(mD\pm E)\cdot C < 0$. 

For $C \in \mathscr{I}(X) \setminus S$, we need to ensure that $(mD\pm
E) \cdot C > 0$.

There  exist only finitely many $C\in \big{(}\mathscr{I}(X)\setminus
S\big{)}$ $\cap$ Neg($m D\pm E$) such that $(m D\pm E)\cdot C<0$. 
Then we can take large enough $n$ such that $((mD\pm E)+nD)\cdot C>0$
for all such $C$,
as $D\cdot C >0$. 

This proves (1).

(2) Clearly, $W_{S'} \cap W_S=\emptyset$ implies $S'\neq S$. 
 
 For the converse, suppose that $S'\neq S$. We may assume, without
 loss of generality, that $ S'\setminus S\neq \emptyset$. Let $C\in
 S'\setminus S$. 
Then $C\cdot D'<0$, $\forall D\in W_{S'}$. Hence $D\notin W_S$. This
proves (2). 
\end{proof}
\end{lem}

Now we show that $W_S$ are precisely the Weyl chambers. 

\begin{lem}\label{lem2.3}
 Connected components of $\text{Big}(X) \setminus \bigcup\limits_{C\in\mathscr{I}(X)}
 C^{\bot}$ are exactly $W_S$ .
 
\begin{proof}
 Let $\text{Big}(X) \setminus \bigcup_{C\in\mathscr{I}(X)} C^{\bot}=\cup
 X_i$, where $X_i$ are the connected components. Let us consider a big
 divisor $D$ such that $D\in X_i$. We define the set
 $S_D=\{C\in\mathscr{I}(X)\vert D\cdot C<0\}$. Note that $S_D$ can be
the empty set. Since $S_D \subset \text{Neg}(D)$, we have $S_D \in
 \mathscr{Z}(X)$. 
 
We claim that $S_D=S_{D'}$ for any two $D, D'\in X_i$. Assume that the
claim is not true. Then with out loss of generality, we can assume 
that there exists $C\in S_{D'}$, but $C\notin S_{D}$, i.e., $C\cdot D>0$ but $C\cdot D'<0$. 
 
 We define a  map $f_C : N^1(X)_{\mathbb{R}}\rightarrow \mathbb{R}$ by
 $f_C(L)=C\cdot L$. It is easy to see that $f_C$ is continuous. 
So $f_C(X_i)$ is a connected subset of $\mathbb{R}$. Note that
$f_C(D)>0$, and $f_C(D')<0$. Thus there exists $D''\in X_i$ such that
$f_C(D'')=0$. 
This contradicts the fact that $X_i$ is in the complement of the set
$\bigcup_{C\in\mathscr{I}(X)}C^\bot$. 
Hence, our assumption is wrong and the claim is proved. So $X_i= W_S$,
where $S=S_D$, for any $D\in X_i$.
\end{proof}
\end{lem}

As a consequence of the above two lemmas, we prove the following
theorem. This result generalizes \cite[Theorem 1.3]{BF} which proves
the same statement for K3 surfaces. 

\begin{thm}\label{thm2.4}
 Let $X$ be a smooth projective variety. Then there is a bijection between the set of Zariski chambers in Big($X$) and the set of Weyl chambers
in Big($X$).
\begin{proof}
 From Lemma \ref{lem2.1}, it follows that there is a bijective
 correspondence between the set of Zariski chambers and $\mathscr{Z}(X)\cup
 \{\emptyset\}$. 
From Lemma \ref{lem2.3}, it follows that there is a bijection 
between set of Weyl chambers and $\mathscr{Z}(X)\cup \{\emptyset\}$. 
\end{proof}

\end{thm}

We use the following result (see \cite[Lemma A.1]{BF})  at several places in subsequent sections.
\begin{lem}\label{lem2.4}
Let $S=(s_{i,j})$ be a negative definite matrix over $\reals$ such that
$s_{i,j} \ge 0$ for $i \neq j$. Then all the entries of the inverse
matrix $S^{ -1}$ are less than or equal to 0.
\end{lem}

\section{Main results}\label{main-results}

We now give necessary and sufficient conditions for a Weyl chamber
to be contained in the interior of a Zariski chamber and vice
versa. These results have been proved for K3 surfaces in \cite{BF}. We
generalize them to an arbitrary surface. Our proofs are similar to the
proofs given in \cite{BF}. 

\begin{thm}\label{prop3.1}
 Let $X$ be a nonsingular projective surface and let $S\in \mathscr{Z}(X)$. We have $W_{S}\subseteq Z_S$ if and only if the following condition holds:
 
 If $C'$ is a curve such that $C'\in \mathscr{I}(X)\setminus S$ and $S\cup \{C'\}\in \mathscr{Z}(X)$, then $C'\cdot C=0$ $\forall C\in S$.
 
 \begin{proof}
  Suppose first that the given condition holds. We will show that
  $W_{S}\subseteq Z_S$. It suffices to show that Neg$(D) =S$ for any
  $D\in W_S$. 
  
 Let $S=\{C_1,\ldots ,C_r\}$ and $D\in W_S$. By the definition of
 $W_S$, we have $D\cdot C_i<0$ for 
$i=1,2,\ldots,r$. Let $D=P_D+N_D$ be the Zariski decomposition of
$D$. As $P_D$ is the nef divisor and 
$C_i$'s are irreducible curves, $P_D\cdot C_i\geq 0$. So $N_D\cdot
C_i<0$ $\forall i$. This implies 
each $C_i$ is a component of $N_D$. Hence $S\subseteq$ Neg($D$).

We will now prove that Neg$(D) \subseteq S$. Suppose not. 
Then Neg($D$)=$\cup_{i=1}^r C_i\cup_{j=1}^s C_j'$ with  $C_j'\notin
S$, $\forall j$ and $s > 0$.  Let $N_D =  \sum\limits_{i=1}^ra_iC_i +
\sum\limits_{j=1}^sa_j'C_j'$.

Then by the given condition, we have $C_i\cdot C_j'=0 ,\forall 1 \le i
\le r, 1 \le j \le s$. As $D\in W_S$ and $C_j'\notin S$, we have $D\cdot C_j'>
0$. 
Let $b_j : = D\cdot C_j'=N_D\cdot C_j' > 0$.  This implies that 
$\sum\limits_{i=1}^sa_i'C_i'\cdot C_j'=b_j$ for $1 \le j \le s$. 

Re-writing the above $s$ equations in matrix form, we get 
$(a_1',\cdots,a_s') A = (b_1,\cdots,b_s)$, where $A$ is the
intersection matrix of $\{C_1',\ldots, C_s'\} \in \mathscr{Z}(X)$.
This gives $(a_1',\cdots,a_s')  =  (b_1,\cdots,b_s) A^{-1}$. Note that
$A$ is a negative-definite matrix, and hence invertible. Moreover, by 
Lemma \ref{lem2.4}, all the entries of $A^{-1}$ are non-positive. 
It follows that $a_j'<0$ for all $j$. But this contradicts the fact
that  $N_D$ is an effective divisor and $C_i'$ are components of $N_D$.  
Thus Neg$(D)=S$ and hence $W_S\subseteq Z_S$.

To prove the other direction of the theorem, suppose that $S =
\{C_1,\ldots, C_r \} \in
\mathscr{Z}(X)$ does not satisfy the given condition. Let $C' \in
\mathscr{Z}(X) \setminus S$ be such that $S \cup \{C'\} \in
\mathscr{Z}(X)$ and $C' \cdot C > 0$ for some $C \in S$. 
We will show that $W_S\nsubseteq Z_S$, by 
exhibiting a divisor $D\in W_S$ which is not in $Z_S$.

The intersection matrix of the set $S\cup \{C'\}$ is negative
definite. 
By Lemma \ref{lem2.1}, there exists a big divisor $D''$ such that Neg$(D'')
= S\cup\{C'\}$.  In fact, there is a bijective correspondence between
the Zariski chambers (different from the nef chamber)  and the elements of
$\mathscr{Z}(X)$. So we may choose a big divisor $D''$ in the interior
of the Zariski chamber corresponding to $S\cup \{C'\}$. By
\cite[Proposition 1.8]{BKS}, we have $\text{Null}(P_{D''}) =
\text{Neg}(D'') =  S\cup\{C'\}$.

We will first construct a big divisor $D' = P_{D'}+ c'C'+\sum\limits_{i=1}^rc_iC_i$,
where $N_{D'}=c'C'+\sum\limits_{i=1}^rc_iC_i$ and $c', c_1, \ldots,
c_r$ are positive real numbers, such that 
$D'\cdot C'=N_{D'}\cdot C' < 0$ and $D'\cdot C_i=N_{D'}\cdot C_i < 0$,
for all $i=1,2,\cdots, r$.

Let  $b', b_1,\cdots ,b_r$ be any negative real numbers. 
Consider the following system of $r+1$ linear equations in $r+1$
variables. 
\begin{eqnarray*}
N_{D'}\cdot C'&=&b'<0, \\
N_{D'}\cdot C_i&=&b_i<0, i=1,2,\cdots,r.
\end{eqnarray*}

Since the intersection matrix of $N_{D'}$ is negative definite, there
exists a unique solution $(c',c_1,\ldots, c_r)$ for the above
system. Moreover, by Lemma \ref{lem2.4},  $c', c_1,\cdots, c_r$ are
positive real numbers. So we have a big divisor $D'$ with the desired properties. Note that
$P_{D'} = P_{D''}$. 

Now we consider the divisor $$D=P_{D'}+\frac{min\{c',c_i\}}{2\vert
  C'^2\vert}C'+\sum\limits_{i=1}^rc_iC_i.$$ 
Then $D$ is a big divisor and its Zariski decomposition is given by
this above decomposition. We claim $D\in W_S$, but $D\notin Z_S$. 

Since Neg$(D)\neq S$, clearly $D\notin Z_S$.

To prove that $D \in W_S$, we compute the intersection number $D \cdot C$
for all negative curves $C \in \mathscr{I}(X)$. 

First, let $C \in S$. Then 
$D\cdot C=\frac{min\{c',c_i\}}{2\vert C'^2\vert}C'\cdot C
+\sum\limits_{i=1}^rc_iC_i\cdot C \le  N_{D'}\cdot C <0$. 

We have 
$D\cdot C'=\frac{min\{c',c_i\}}{2\vert
  C'^2\vert}C'^2+\sum\limits_{i=1}^rc_iC_i\cdot
C'=\frac{-\text{min}\{c',c_i\}}{2}+\sum\limits_{i=1}^rc_iC_i\cdot
C'>0$, as  $C'\cdot C_i>0$ for some $i$.


Now let $C \notin  S\cup \{C'\}$. Then $P_D \cdot C > 0$, since $D''$
is the interior of $Z_{S\cup \{C'\}}$ and $P_D = P_{D''}$. So $D \cdot
  C = P_D \cdot C + N_D \cdot C > 0$. 

This proves that $D \in W_S$ and completes the proof the theorem. 
\end{proof}
\end{thm}


We will now give a criterion for the interior of a Zariski chamber to
be contained in a Weyl chamber. 

\begin{thm}\label{prop3.2}
Let $X$ be a smooth projective surface and let $S\in \mathscr{Z}(X)$. Let $\mathring{Z}_S$ be the interior of the Zariski chamber $Z_S$. Then we have  $\mathring{Z}_S\subset W_S$ if and only if $C'\cdot C''=0$ for all curves $C', C'' \in S$.

\begin{proof}
First, consider a set $S =\{ C_1,\cdots, C_r\} \in \mathscr{Z}(X)$ satisfying the given condition. We will show
that $\mathring{Z}_S\subset W_S$. Let $D \in \mathring{Z}_S$.

Let the Zariski decomposition of $D$ be $D=P_D+N_D$, where
$N_D=\sum\limits_{i=1}^r a_i C_i$.  
Then $D\cdot C_i= N_D\cdot C_i= a_i C_i^2<0$, for all
$i=1,\cdots,r$. 
Let $C$ be a negative curve on $X$ which is not in $S$. 
Then $D\cdot C\geq 0$. If $D\cdot C = 0$, then $D \cdot C = P_D\cdot
C+ N_D\cdot C=0$, 
which implies that $P_D\cdot C = N_D\cdot C=0$. Hence, $C\in$
Null$(P_D)$. This shows that  Neg$(D) \ne  \text{Null}(P_D)$. But
since  $D$ is in the interior of the Zariski
Chamber $Z_S$, we have Neg($D$)=Null($P_D$), by\cite[Proposition 1.8]{BKS}. 
Hence we must have $D \cdot C > 0$. So $\mathring{Z}_S\subset W_S$.
 
To prove the other direction of the theorem, consider $S \in
\mathscr{Z}(X)$ such that there exist curves $C_i,C_j\in S$ satisfying
$C_i\cdot C_j\neq 0$. We will show that there exists a big 
divisor $D\in \mathring{Z}_S$ such that $D\notin W_S$. Let $S=\{
C_1,...,C_r\}$ and suppose that $C_1\cdot C_2\neq 0$. 
 
Fix an ample divisor $H$.  For some unknown positive real numbers
$a_1,\cdots,a_r, a_1^{\star}, \cdots, a_r^{\star}$, consider the
following divisor:
 \begin{align*}
  D'=H+ \sum\limits_{i=1}^r a_i^{\star}C_i+\sum\limits_{i=1}^r(a_i-a_i^{\star})C_i,
 \end{align*}
 where $P_{D'}=H+ \sum\limits_{i=1}^r a_i^{\star}C_i$ and $N_{D'}=\sum\limits_{i=1}^r(a_i-a_i^{\star})C_i$.
 Note that $D'$ is a big divisor, being the sum of an ample divisor and
 an effective divisor.

 Now we want to find values of $a_i, a_i^{\star}$, such that
 Neg$(D')=$Null$(P_{D'})$. To obtain this, we solve
 the following system of $r$ linear equations in $r$ variables: 
\begin{align*}
 0=P_{D'}\cdot C_j= H\cdot C_j+\sum\limits_{i=1}^r a_i^{\star} C_i.C_j.
\end{align*}
This system has a unique solution $(a_1^{\star}, \cdots, a_r^{\star})$
with $a_i^{\star} > 0$ for all $1 \le i \le r$. 
Choose real numbers $a_i$ such that $0< a_i^{\star} <a_i$ for all
$i$. Finally, choose a real number $k$ such that $0 <
k<\frac{(a_2-a_2^{\star})}{\mid C_1^2\mid }$ and consider 
the divisor:
\begin{align*}
 D=H+(a_1^{\star}+k)C_1+\sum\limits_{i=2}^r a_i C_i.
\end{align*}
Then the Zariski decomposition of $D$ is given by $D=
P_D + N_D$, 
where $P_{D}=P_{D'}=H+ \sum\limits_{i=1}^r a_i^{\star}C_i$ and $N_D=k
C_1+\sum\limits_{i=2}^r(a_i-a_i^{\star})C_i$. 
As Neg($D$) = Null($P_{D}$) = $S$, we have $D\in \mathring{Z}_S$.
But $D \cdot C_1=k C_1^2+\sum\limits_{i=2}^r(a_i-a_i^{\star})C_i. C_1\geq k
C_1^2+ (a_2-a_2^{\star})>0$, by the choice of $k$. So $D \notin W_S$. 

This completes the proof. 
\end{proof}
\end{thm}


We recall the following result in \cite{RS} which determines when
the interior of each Zariski chamber coincides with a Weyl chamber.

\begin{thm}\cite[Theorem 3]{RS} \label{prop3.3}
Let $X$ be a smooth projective surface. The following conditions are equivalent:

(a) the interior of each Zariski chamber on $X$ coincides with a simple Weyl chamber,

(b) if two different irreducible negative curves $C_1 \neq C_2$ on $X$ meet (i.e., $C_1 \cdot C_2 > 0$), then
\begin{align*}
 C_1 \cdot C_2 \geq \sqrt{C_1^2 \cdot C_2^2} .
\end{align*}
\end{thm}

If the above equivalent conditions hold on a surface $X$, then
\cite{RS} says that the 
Zariski chambers of $X$ are \textit{numerically
  determined}. 

Theorems \ref{prop3.1} and \ref{prop3.2}, determine
when a specific Weyl chamber is contained in a Zariski chamber and
when the interior of a specific Zariski chamber is contained in a Weyl
chamber. Our results imply Theorem \ref{prop3.3}, as we show
below. 

We first give some equivalent formulations of condition (b)
of Theorem \ref{prop3.3}. 
 
 \begin{thm}\label{thm3.6} Let $X$ be a smooth projective surface. 
  The following are equivalent.
\begin{enumerate} 
\item If two irreducible negative curves $C_1 \neq C_2$ on $X$ meet (i.e., $C_1 \cdot C_2 > 0$), then
$ C_1 \cdot C_2 \geq \sqrt{(C_1^2 \cdot C_2^2)}$. 
\item If $C_1$ and $C_2$ are any two negative curves such that $\{
  C_1, C_2\}\in \mathscr{Z}(X)$, then $C_1\cdot C_2=0$.
\item Let $S \in \mathscr{Z}(X)$. If $C'\in \mathscr{I}(X)\setminus
  S$, and $S\cup C'\in \mathscr{Z}(X)$, then $C'\cdot C=0$ for all
  curves $C\in S$.
\item Let $S \in \mathscr{Z}(X)$. Then $C_1\cdot C_2=0$ for all curves
  $C_1, C_2 \in S$. 
\end{enumerate}

\begin{proof}

\underline{(1)$\Rightarrow$ (2)}:
  Assume that (2) does not hold. Let  $S=\{ C_1, C_2\}\in
  \mathscr{Z}(X)$ but $C_1\cdot C_2 > 0$. 
As the intersection matrix of $S=\{ C_1, C_2\}$ is negative definite, we have,
\begin{align*}
 x^2 C_1^2+ y^2C_2^2+2xy C_1\cdot C_2 <0
\end{align*}
for any $(x,y) \ne (0,0)$.
For any $(x,y)\neq (0,0)$, we have
\begin{align*}
 C_1\cdot C_2<\frac{(x\sqrt{-C_1^2})^2+(y\sqrt{-C_2^2})^2}{2xy\sqrt{C_1^2 C_2^2}}\sqrt{C_1^2 C_2^2}
\end{align*}

Take $x=\sqrt{-C_2^2}$ and $y=\sqrt{-C_1^2}$. Then
$C_1C_2<\sqrt{C_1^2C_2^2}$. This violates (1). 

\underline{(2)$\Rightarrow$ (1):}
Let $C_1, C_2$ be two different negative curves. If the intersection
matrix of $S=\{C_1, C_2\}$ is negative definite, then $C_1\cdot
C_2=0$, by (2). So (1) clearly holds. If the intersection matrix of
$S$ is not negative definite, 
then there exists a tuple $(x,y)$, where $x,y$ are both nonzero real
numbers having the same sign, such that
\begin{align*}
 x^2 C_1^2+ y^2C_2^2+2xy C_1\cdot C_2 \geq0
\end{align*}
This implies,
\begin{align*}
 C_1\cdot C_2\geq \frac{(x\sqrt{-C_1^2})^2+(y\sqrt{-C_2^2})^2}{2xy\sqrt{C_1^2 C_2^2}}\sqrt{C_1^2 C_2^2}
\end{align*}

If $a,b$ are positive integers, then $(a-b)^2\geq 0\Rightarrow
\frac{a^2+b^2}{2ab}\geq 1$. So 
the above inequality implies  $C_1 \cdot C_2 \geq \sqrt{C_1^2 \cdot C_2^2}$.

\underline{(3)$\Rightarrow$ (2) and {(4)$\Rightarrow$ (2)}:} These
implications are straightforward.

\underline{(2)$\Rightarrow$ (3)}: Let $S \in \mathscr{Z}(X)$. 
If $S\cup C'\in \mathscr{Z}(X)$ for some $C'\in
\mathscr{I}(X)\setminus S$, and $C \in S$, then $\{C, C' \} \in
\mathscr{Z}(X)$. So it follows that $C \cdot C' = 0$, as required.

\underline{(2)$\Rightarrow$ (4)}: Let $S \in \mathscr{Z}(X)$. If $C_1,
C_2 \in S$, then $\{C_1, C_2 \} \in \mathscr{Z}(X)$. So it follows that $C_1 \cdot
C_2 = 0$, as required. 
\end{proof}
 \end{thm}

\begin{xrem}
\rm 
We will now explain why our main results (Theorems \ref{prop3.1} and \ref{prop3.2}) imply Theorem
\ref{prop3.3}.

First, suppose that the condition (b) of Theorem
\ref{prop3.3} holds. This is same as the statement (1) of Theorem \ref{thm3.6}.
So both statements (3) and (4) of Theorem \ref{thm3.6} hold and
this implies that 
the condition of Theorems \ref{prop3.1}  and \ref{prop3.2} hold for all subsets $S \in \mathscr{Z}(X)$. So for
every $S \in \mathscr{Z}(X)$, we have the equality $W_S =
\mathring{Z}_S$. 

Conversely, if the condition (a) of Theorem \ref{prop3.3} holds, then we have 
 $W_S = \mathring{Z}_S$ for every $S \in \mathscr{Z}(X)$. By Theorem
 \ref{prop3.1}, we conclude that the statement (3) of Theorem
 \ref{thm3.6} is true and hence condition (b) of Theorem \ref{prop3.3}
 follows. 

We will give some examples later which show that Theorem \ref{prop3.3} does not imply
Theorem
\ref{prop3.1} or Theorem \ref{prop3.2} (see Example \ref{exm2}, for instance).
\end{xrem}

Now we give a necessary and sufficient condition for a Weyl chamber
and a Zariski chamber to have non-empty intersection.

\begin{thm}\label{thm3.7}
 Let $X$ be a smooth projective surface and let $S, S_1 \in \mathscr{Z}(X)$. Then
 \begin{align*}
  W_{S_1}\cap Z_S\neq \emptyset
 \end{align*}
if and only if $S_1\subseteq S$ and any subset $S'\subseteq S\setminus
S_1$ satisfies the following property:
there exist $C'\in S'$ and $C \in S\setminus S'$ such
that $C'\cdot C>0$. 
\begin{proof}
 First assume that $ W_{S_1}\cap Z_S\neq \emptyset$.  Let $D\in  W_{S_1}\cap
 Z_S$ and let its Zariski decomposition be $D=P_D+N_D$. Then 
for every curve $C'\in S_1$, we have $D\cdot C'<0$, as $D\in
W_{S_1}$. Moreover, $N_D\cdot C'<0$, as $P_D$ is nef. This implies
that $C'$ is an irreducible component of $N_D$. Since the set of irreducible
components of $N_D$ is precisely $S$ (since $D \in Z_S$), we conclude $S_1\subseteq S$.

Now suppose that there exists a subset $S'\subseteq S\setminus S_1$
 which does not 
satisfy the given condition. In other words, for any $C'\in S'$, we
have $C'\cdot C=0$, for all $C\in S\setminus S'$. The Zariski decomposition of
$D$ can be written as 
$$D=P_D+N_D=P_D+\sum_{C_i\in S'} a_iC_i+\sum_{C_i\in S\setminus
  S'}b_iC_i.$$ 

As $D\in W_{S_1}$, we have $D\cdot C_j>0$ for all $C_j\in S'\subseteq
S\setminus S_1$. Note that,
$$D\cdot C_j=\sum_{C_i\in S'} a_iC_i\cdot C_j> 0.$$

Since $S' \subset S \in \mathscr{Z}(X)$, the intersection matrix of
$S'$ is negative definite. Then it follows from Lemma \ref{lem2.4}
that $a_i<0$, 
which is absurd since $C_j\in S'$ are irreducible component of the
divisor $D\in Z_S$. 
Hence, our assumption on $S'$ is wrong and this completes the proof of 
one direction of the theorem.

 To prove the converse direction, let $S' \subseteq S$ be a subset
 satisfying the given condition. 
Our goal is to find a $D\in  W_{S_1}\cap Z_S$.

First we show that $W_{S_1}\cap Z_{S_1} \ne \emptyset$. We will construct
a divisor in this intersection by defining a valid Zariski
decomposition.  Fix an ample divisor $H$.  
For some positive real numbers $a_i$ and $a_i^*$ (to be determined), consider a
divisor of the following form:
\begin{equation}\label{d1}
D_1=H+ \sum_{C_i\in S_1} a_i^*C_i + \sum_{C_i\in S_1} (a_i-a_i^*)C_i,
\end{equation}
where $P_{D_1} :=H+\sum_{C_i\in S_1} a_i^*C_i$ and $N_{D_1} :=
\sum_{C_i\in S_1} (a_i-a_i^*)C_i$. 
We will show that  $a_i^*$ and $a_i$ can be chosen such that $D_1\in W_{S_1}\cap Z_{S_1}$.
 
If the above decomposition is to be the Zariski decomposition of $D_1$, then we must have
$P_{D_1}\cdot C_j=0$, 
for all $C_j\in S_1$. Hence,
 \begin{align*}
 \sum_{C_i\in S_1} a_i^*C_i\cdot C_j=-H\cdot C_j<0.
 \end{align*}
As the intersection matrix of $S_1$ is negative definite, we can solve
uniquely for $a_i^*$ in the above linear system of equations. Moreover, by 
Lemma \ref{lem2.4}, $a_i^*>0$. In particular, we conclude that $P_D$ is nef. 

Fix a set $\{y_i<0\mid 1\leq i\leq \vert S_1\vert\}$ of negative
real numbers and consider the following linear system of $\vert S_1 \vert$ equations in
$\vert S_1 \vert$ variables,
\begin{equation}\label{ws}
 \sum_{C_i\in S_1} x_i C_i \cdot C_j=y_j, \text{~where~} C_j\in S_1
\end{equation}
Again by the negative definiteness of the intersection matrix of
$S_1$ and Lemma \ref{lem2.4}, there is a unique solution  $x_i>0$ for
the above system.  
Set $a_i := a_i^*+x_i$. Now we conclude that $D_1\in W_{S_1}\cap
Z_{S_1}$, where $D_1$ is defined in \eqref{d1}. Indeed,  $D_1$ is big,
since it is the sum of an ample divisor and an effective
divisor. Further, $D_1 =
P_{D_1}+N_{D_1}$ is a Zariski decomposition by construction. Hence
$\text{Neg}(D_1) = S_1$, so $D_1 \in Z_{S_1}$. Also, by \eqref{ws},
$D_1 \cdot C_j < 0$ for every $S_j \in S_1$. If $C \notin S_1$, then 
$D_1 \cdot C \ge H \cdot C > 0$.  So $D_1 \in W_{S_1}$, as required.

Now, define 
\begin{center}
 $S_2 : =\{C_2\in S\setminus S_1\mid C_2\cdot C_1>0$ for some $C_1\in S_1\}$.
\end{center}
By the condition given in the theorem, $S_2\neq \emptyset$. By following
a similar construction as above,  we find a divisor $D_2\in W_{S_1}\cap Z_{S_1\cup S_2}$.
We briefly explain the procedure, for clarity. Let 
\begin{center}
 $D_2= H+\sum_{C_i\in S_1}\alpha_i C_i+\sum_{C_i\in S_2} \beta_iC_i+n N_{D_1}+ \sum_{C_i\in S_2} C_i$,
\end{center}
where $P_{D_2}=H+\sum_{C_i\in S_1}\alpha_i C_i+\sum_{C_i\in S_2}
\beta_iC_i$ and $N_{D_2}=n N_{D_1}+ \sum_{C_i\in S_2} C_i$. We claim
that there exist positive real numbers $\alpha_i, \beta_i$ and a
sufficiently large integer $n$ such that 
$D_2\in W_{S_1}\cap Z_{S_1\cup S_2}$. 

We can set-up linear equations as above to find positive real numbers
$\alpha_i, \beta_i$ such that $D_2 =P_{D_2}+N_{D_2}$ is the Zariski
decomposition of the big divisor $D_2$ for any positive integer $n$. 
It is also clear then $D_2 \in Z_{S_1 \cap S_2}$. 
Now we show that for
a sufficiently large integer $n$, $D_2 \in W_{S_1}$. 

Let $C \in S_1$. Then $D_2 \cdot C = nN_{D_1}\cdot C + \sum_{C_i\in
  S_2} C_i\cdot C$. Since $N_{D_1} \cdot C < 0$, we may choose $n \gg
0$ such that $D_2 \cdot C < 0$ for every curve $C \in S_1$. Now
let $C \notin S_1$. If $C \notin S_2$, then clearly $D_2 \cdot C >
0$ (note that $D_2$ is a sum of the ample divisor $H$ and an effective
divisor supported on the curves in $S_1 \cap S_2$). If $C \in S_2$,
then, by the definition of $S_2$,  there exists $C' \in S_1$ such
that $C \cdot C' > 0$. Now $D_2 \cdot C = n N_{D_1} \cdot C +
\sum_{C_i\in S_2} C_i \cdot C$. Since all the terms in $N_{D_1} \cdot
C$ are non-negative and at least one term is positive, we can choose
$n \gg 0$ such that $D_2 \cdot C > 0$. 

Proceeding this way, since $S$ is a finite set, we find a divisor
$D\in W_{S_1}\cap Z_S$, as claimed. 
This completes the proof of the theorem. 
\end{proof}

\end{thm}

\subsection{Examples  and remarks}
We now give some examples and make some remarks illustrating our results.

\begin{exm}\rm
Let $\pi: X \to \mathbb{P}^2$ be the blow up at four collinear points $P_1,
P_2, P_3, P_4\in \proj^2$. Let $H$ denote the pull-back of
$\str_{\proj^2}(1)$ and let $E_i = \pi^{-1}(p_i)$ be the exceptional
divisors.

The set of irreducible negative curves of
$X$ is given by: 
$$\mathscr{I}(X) = \{E_1, E_{2}, E_{3}, E_{4},\tilde{L}_{1234}\},$$
where $\tilde{L}_{1234}$ is the strict transformation of the line joining
$P_1, P_2, P_3,$ and $P_4$. Note that $\tilde{L}_{1234} =
H-E_1-E_2-E_3-E_4$. 

Let $S = \{E_1, E_2, \tilde{L}_{1234}\}$ and $S_1 = \{E_1, E_2\}$. It
is easy to check that the intersection matrices of both $S$ and $S_1$
are negative definite. So $S, S_1 \in \mathscr{Z}(X)$. Using Theorem
\ref{thm3.7}, we can conclude that $W_{S_1}\cap Z_S \neq \emptyset$. We
exhibit below an explicit divisor in the intersection. 

Consider the divisor $D = 6H - E_{3}- E_{4} + E_{1} + E_{2}$. 
 Note that the Zariski decomposition of $D$ is given by
$$D = (5H) + (\tilde{L}_{1234} + 2E_{1} + 2E_{2}).$$

It follows from the construction that $D\in Z_S$ (note that $D
\notin \mathring{Z}_S$). Note that, $D \cdot
E_{1} = -1$, $D \cdot E_{2} = -1$, while $D \cdot E_{3} = D \cdot
E_{4} = 1$ and $D\cdot \tilde{L}_{1234} = 6$. Hence $D \in W_{S_1}$ too.
\end{exm}

\begin{exm}\label{exm2}
\rm
 Let $X \to \mathbb{P}^2$ be a blow up at five points $P_1,P_2,P_3,
 P_4, P_5 \in \proj^2$ such that $P_1,P_2,P_3$ are collinear and no
 other triple is collinear. Let $C: = \tilde{L}_{P_1P_2P_3}$ denote the 
strict transformation of the line containing $P_1$, $P_2,
P_3$. Then $C^2 = -2$ and $C \cdot E_1 = 1$. So we see that the
condition of Theorem \ref{prop3.3} is not satisfied. In other words,
$X$ is not numerically determined. We exhibit two sets in
$\mathscr{Z}(X)$ which behave differently with regard to the
containment of Weyl
and Zariski chambers. 

It is easy to verify that the set $\mathscr{I}(X)$ of negative curves
on $X$ consists of the exceptional divisors, $C$ and lines through all
pairs $P_i, P_j$ of points, where either $P_i$ or $P_j$ is not in
$\{P_1,P_2,P_3\}$. 

Let $S=\{E_1, \tilde{L}_{P_1P_2P_3}\}$ and $S'=\{E_{4},
E_{5}\}$. Note that $S, S'\in \mathscr{Z}(X)$.

It is clear that $S'$ satisfies the condition of Theorem \ref{prop3.2}. It can
be checked that $S'$ also satisfies the condition of Theorem
\ref{prop3.1}. If $C'$ is a negative curve which meets either $E_4$ or
$E_5$, then it turns out that $\{E_4,E_5,C'\} \notin
\mathscr{Z}(X)$. For example, if $C' = H-E_4-E_5$, then the
intersection matrix of $\{E_4,E_5,C'\} $ has a positive eigenvaule. On
the other hand, the intersection matrix of $\{E_4,E_5,H-E_1-E_4\}$ is
not even invertible. Thus we have $\mathring{Z}_{S'}=W_{S'}$. 

On the other hand, the condition of Theorem \ref{prop3.2} fails for
$S$. So $\mathring{Z}_{S}\not \subset W_{S}$. And we can check that 
the condition of Theorem \ref{prop3.1} holds for $S$. So 
$\mathring{Z}_{S}\supsetneq W_{S}$. In fact, by Theorem \ref{thm3.7},
$Z_S \cap W_{S''}  = \emptyset$ for any $S'' \ne S$. 
 \end{exm}

\begin{exm}\label{exm3}
\rm Let $X \to \proj^2$ be a blow up at ten points of intersection of
five general lines $L_1, \cdots, L_5$ in $\proj^2$. Suppose that the
four points of intersection that lie on $L_1$ are $p_1,\ldots,
p_4$. Then the strict transform $C_1 : = H-E_1-E_2-E_3-E_4$ of $L_1$ is a
negative curve. Let $C_2 : = E_1$. Then $S= \{C_1,C_2\} \in
\mathscr{Z}(X)$. Since $C_1\cdot C_2 = 1$ and $C_1^2 = -3$, the
condition in Theorem \ref{prop3.3} is not satisfied. So again the
Zariski chambers are not numerically determined. 

In fact, using Theorems \ref{prop3.1}  and  \ref{prop3.2}, we see that
neither of the inclusions $W_S \subset Z_S$ nor
$\mathring{Z}_{S}\subset W_{S}$ hold. Indeed, since $C_1 \cdot C_2 \ne 0$,
we know by Theorem \ref{prop3.2} that $\mathring{Z}_S \not\subset W_S$. On other hand, if $C = E_2$,
then it is easy to check that $S \cup \{C\} \in \mathscr{Z}(X)$. Since 
$C \cdot C_1 = 1 \ne 0$, we know by Theorem \ref{prop3.1} that 
$W_S \not\subset Z_S$. 

On the other hand, by Theorem \ref{thm3.7}, we know that $W_S \cap Z_S
\ne \emptyset$.

\end{exm}

\begin{exm}\label{exm4}\rm

Let $D$ be an irreducible and reduced plane cubic and let $X \to
\proj^2$ be a blow up of $s$ very general points on $D$. 
It is well-known that the only 
the only negative curves on $X$ are the strict
transform $C$ of $D$ (when $s > 9$), and the \textit{$(-1)$-curves},
i.e., 
smooth rational curves whose self-intersection is $-1$; see
\cite{H} or \cite{D-et-al} (refer to the arXiv version of the latter paper for
substantial changes made after publication). 
Since $C^2 = 9-s$ and $C \cdot C' = 1$ for any other negative curve
$C'$, the Zariski chambers on $X$ are not numerically determined for
$s \ge 11$. 

In fact, it is easy to check that $S = \{C, C'\} \in \mathscr{Z}(X)$
for any negative curve $C' \ne C$. Since $C \cdot C' =1$, we have
$\mathring{Z}_S \not\subset W_S$, by Theorem \ref{prop3.2}. On the other hand, $S \cup \{C''\}
\notin \mathscr{Z}(X)$ for any negative curve $C'' \notin S$. So 
$W_S \subset Z_S$, by Theorem \ref{prop3.1}. 
\end{exm}

 \begin{xrem}\label{delpezzo}
  \rm Let $X$ be a Del Pezzo surface. Then $X$ is $\proj^2$, $\proj^1
  \times \proj^1$, or a blow up of $\proj^2$ at eight or fewer general
  points. The number of Zariski chambers in each of these cases is
  known; see \cite[Theorem]{BFN}. 

There are no negative curves in $\proj^2$ or $\proj^1
  \times \proj^1$.  If $X$ is a blow up of $\proj^2$ at eight or fewer general
  points, it is well-known that the only negative curves on $X$ are
  $(-1)$-curves.  Hence the condition given in \cite[Theorem
  3]{RS} (see Theorem \ref{prop3.3}) is satisfied. So the interior of each Zariski chamber coincides with 
a Weyl chamber.
 \end{xrem}

\begin{xrem}\rm
Let $X \to \proj^2$ be a blow up of $r \ge 0$ general points of
$\proj^2$. The \textit{$(-1)$-curves Conjecture } (some times called
the \textit{Weak SHGH Conjecture}) predicts that the only negative
curves on $X$ are the $(-1)$-curves. This is known to be true when $r
\le 9$ (see Remark \ref{delpezzo} above) but it is open when $r \ge
10$. If this conjecture is true, then the Zariski chambers on $X$ are
numerically determined. 
\end{xrem}

\begin{xrem}
 \rm Let $X$ be a geometrically ruled surface over a nonsingular curve $Y$. The
 Zariski chamber and Weyl chamber decomposition of the big cone of $X$
 coincide. If $\mathcal{E}$ is a semistable bundle over $Y$, then we
 know that the big cone and the ample cone are the same in $X =
 \mathbb{P}_Y(\mathcal{E})$. 
So we have only one Zariski chamber and one Weyl chamber. If
$\mathcal{E}$ is an 
unstable bundle over $Y$, then $X$ has exactly one negative curve. Hence
there are two Zariski chambers and two Weyl chambers. One of the
Zariski chambers is the nef chamber and one of the Weyl chambers is
the ample cone (corresponding to $S= \emptyset$). The other chamber corresponds to the unique negative
curve on $X$. 
\end{xrem}

\begin{xrem}\rm
The Weyl and Zariski chamber decomposition on K3 surfaces is studied
in detail in \cite{BF}. Our results (Theorems \ref{prop3.1} and
\ref{prop3.2}) and their proofs are motivated by analogous results
proved for K3 surfaces. In addition, \cite{BF} gives examples of K3
surfaces where the decompositions coincide and where the decomposition
differ. 

The case of Enriques surfaces is considered in \cite{RS}. The authors
relate the coincidence of Weyl chambers and interiors of Zariski
chambers to the properties of elliptic fibrations on an Enriques
surface.
\end{xrem}

\end{document}